\documentclass[11pt]{article}
\usepackage{ulem}
\usepackage[makeroom]{cancel}
\usepackage{ulem}
\usepackage[usenames,dvipsnames,svgnames,table]{xcolor}
\usepackage{xparse}
\usepackage{tikz}
\usetikzlibrary{positioning,patterns,matrix,arrows,calc}
\usetikzlibrary{decorations.markings}
\usetikzlibrary{patterns}
\usetikzlibrary{snakes}
\usetikzlibrary{shapes}
\usepackage[latin1]{inputenc}
\usetikzlibrary{trees}
\usepackage{hyperref} 
\usetikzlibrary{decorations.pathmorphing}
\usepackage{amsfonts, simplewick}
\usepackage{amssymb}
\usepackage{amscd}
\usepackage{amsmath, mathtools}
\usepackage{epsfig}
\usepackage{graphicx, subfigure}
\usepackage{latexsym}
\usepackage[font=small,labelfont=bf]{caption} 
\usepackage{float}
\usepackage{verbatim}
\usepackage{tipa}
\usepackage{amsthm,amssymb}
\usepackage{color}
\usepackage{diagbox}
\newlength{\dinwidth}
\newlength{\dinmargin}
\newcommand{\be}{\begin{equation}}
\newcommand{\ee}{\end{equation}}
\newcommand{\ba}{\begin{eqnarray}}
\newcommand{\ea}{\end{eqnarray}}

\newtheorem{definition}{Definition}
\newtheorem{theorem}{Theorem}
\newtheorem{conjecture}{Conjecture}
\newtheorem{proposition}{Proposition}
\newtheorem{corollary}{Corollary}
\newtheorem{remark}{Remark}
\newtheorem{lemma}{Lemma}
\newtheorem{example}{Example}

\def\e{{e}}
\begin{document}
\title{Enumeration of multi-rooted plane trees }
\author{Anwar Al Ghabra$^1$, K. Gopala Krishna$^2$,   Patrick Labelle$^3$, and Vasilisa Shramchenko$^1$}

\footnotetext[1]{Department of mathematics, University of
	Sherbrooke, 2500, boul. de l'Universit\'e,  J1K 2R1 Sherbrooke, Quebec, Canada. E-mails: {\tt mouhammed.anwar.al.ghabra@usherbrooke.ca} and {\tt Vasilisa.Shramchenko@Usherbrooke.ca}}

\footnotetext[2]{E-mail: {\tt Gopala.K.Krishna@gmail.com}}

\footnotetext[3]{E-mail: {\tt patrick.labelle3@usherbrooke.ca}}

\date{}
\maketitle

 \begin{abstract}
We give closed form expressions for the numbers of multi-rooted plane trees with specified degrees of root vertices. This results in an infinite number of integer sequences some of which are known to have an alternative interpretation. We also propose recursion relations for numbers of such trees as well as for the corresponding generating functions. Explicit expressions for the generating functions corresponding to plane trees having two and three roots are derived. As a by-product, we obtain  a new binomial identity and  a conjecture  relating hypergeometric functions. 
 \end{abstract}
 
 \vskip 1cm

\\


\tableofcontents

\section{Introduction}

A plane tree is a particular case of a  connected  
 ribbon graph, or a map, that is a graph embedded into a compact orientable surface in such a way that every face is homeomorphic to a disc. The genus of the surface is also called genus of the embedded graph.  A tree is embedded in this way into a sphere.
\\

By assigning lengths to edges of a ribbon graph, one obtains a metric ribbon graph. The spaces of metric ribbon graphs with vertices of degree three or higher give a way to describe combinatorially the moduli spaces of Riemann surfaces with marked points, see \cite{MP}, the bridge between metric ribbon graphs and Riemann surfaces being obtained by the theory of Strebel differentials \cite{Strebel}.  Initially, ribbon graphs were used to describe spaces of Riemann surfaces by R. Penner \cite{Penner}. In a seminal paper \cite{thooft}, 't Hooft showed how, in  a certain limit, Feynman diagrams of non-abelian gauge theories can be analyzed using ribbon graphs. This has also led to connections between point particle quantum field theories and string theory.  Ribbon graphs also arise naturally in the context of matrix models for quantum field theories, see for example \cite{Zvonkin}. 
\\

 Bipartite ribbon graphs, in particular plane trees, can be seen \cite{Grothendieck} as representing Belyi pairs, that is pairs of a Riemann surface and a meromorphic function on this surface with critical values in the set $0,1,\infty.$ Such graphs are called {\it dessins d'enfant} following Grothendieck. In \cite{Ka, Zo} it is shown that the generating function of numbers of dessins d'enfant satisfies the KP (Kadomtsev-Petviashvili) hierarchy and in \cite{KaZo-Virasoro}
the same function is shown to satisfy Virasoro constraints and the {\it topological recursion of Chekhov-Eynard-Orantin} \cite{EO} for an appropriate spectral curve.   Plane trees are included in this generating function as dessins d'enfant representing Belyi pairs given by a Riemann sphere and a polynomial function. On the other hand, any ribbon graph can be seen as representing a so-called {\it clean Belyi pair} \cite{Dumi} by introducing an extra vertex of degree two in the middle of every edge.  Numbers of clean Belyi pairs are linked to the topological recursion in \cite{Dumi}.  Enumeration of ribbon graphs is closely related to computing Hurwitz numbers, see for example \cite{GouldenJackson, DubrovinYangZagier} and \cite{GouldenJackson} for the relationship with the KP hierarchy. 
\\

The question of enumeration of maps was first considered by W. T. Tutte in \cite{Tutte}. In order to simplify the problem of counting, Tutte introduced a {\it root} in a graph, that is a distinguished orientation of one of the edges. Two graphs are identified if they can be obtained from one another by a homeomorphism of the underlying surface in such a way that the distinguished edge is mapped to the distinguished edge and the orientations agree.  The presence of a root thus ensures that the graph has no symmetries, that is no nontrivial automorphisms. 
The following formula for the number of rooted maps on the sphere having $e$ edges was derived in  \cite{Tutte}: 
\begin{equation}
\label{Tutte}
\frac{2(2e)!3^e}{e!(e+2)!}\,.
\end{equation}
This formula combines all graphs of genus zero having $e$ edges, some of which are plane trees. Following the seminal paper of  Tutte, numerous results on enumeration of rooted maps appeared.   Various techniques were developed for such enumeration, including recursive construction of maps and deducing differential equation on the generating functions, see \cite{AB} for a partial overview. 
Here we give a brief review of the results most closely related to the approach of our paper. 
\\

Tutte's result  \eqref{Tutte} was generalized in \cite{AB} to give the number $m_1(e)$ of rooted graphs with $e$ edges, combining maps of all genera, in the form: 
\begin{equation}
\label{ABm1}
m_1(e) = \frac{1}{2^{\e+1}} \sum_{i=0}^{\e} (-1)^i \sum_{\substack{k_1+\cdots + k_{i+1}= \e+1 \\ k_1, \ldots, k_{i+1} > 0}} \,\,\prod_{j=1}^{i+1}  \frac{(2k_j)!}{k_j!} \;.
\end{equation}

The following result of  \cite{WalshLehman-1} can be used to separate the numbers of graphs by genus. The formula of \cite{WalshLehman-1} is given in terms of the numbers $C_{g,v}$, which were called in \cite{Motohico_Catalan} the {\it generalized Catalan numbers}. They are defined as follows. The integer $C_{g,v}(d_1,\dots, d_v)$ is a number of ribbon graphs of genus $g$ with $v$ ordered vertices such that the vertex number $j$ is incident to $d_j$ half-edges; moreover, at each vertex one of the incident half-edges is marked. Such graphs are called {\it dicings} in \cite{WalshLehman-1}. The graph with one vertex and no edges is also considered as a degenerate dicing, giving $C_{0,1}(0)=1\,.$ The Catalan numbers $C_m = \frac{1}{m+1}  {2m\choose m}$ are obtained as a particular case, namely $C_m=C_{0,1}(2m)$. In other words, the $m$th Catalan number is the number of genus zero maps with one vertex and $m$ edges, where one edge is given an orientation. As before, two graphs are identified if they can be mapped to each other by a homeomorphism of the sphere preserving the chosen orientation of the marked edge.   The number $ m_1(e,v;g)$ of rooted maps of genus $g$ with $v$ vertices and $e$ edges is then obtained in the form \cite{WalshLehman-1}:
\begin{equation}
\label{m1Cgn}
 m_1(e,v;g)=\frac{2e}{v!}  \sum_{d_1+d_2+\dots+d_v=2e}  \frac{C_{g,v}(d_1,\dots, d_v)}{d_1\dots d_v}\;.
\end{equation}

The generalized Catalan numbers $C_{g,v}(d_1,\dots, d_v)$ can be obtained using a recurrence relation from \cite{WalshLehman-1} or \cite{Motohico_Catalan}. Alternatively, they can be computed by the algorithm of the topological recursion of Chekhov-Eynard-Orantin \cite{EO} applied to the algebraic curve corresponding to the equation $y^2=x^2-2$, see \cite{Motohico_Catalan} and also \cite{KriPatVas2}. 
\\

As mentioned, introducing one distinguished half-edge, a root,  removes the possibility of non-trivial automorphisms of a map and simplifies enumeration. Introducing further roots is thus unnecessary from the point of view of destroying symmetry. However, it is also interesting to consider multi-rooted maps. For example, dicings from \cite{WalshLehman-1} are multi-rooted maps as each their  vertex is a {\it root-vertex}, that is incident to a root. We will refer to such maps as maximally rooted. Furthermore, in    \cite{AG, AG2, BCR, KriPatVas1}, $N$-rooted maps were considered, that is maps having $N$ ordered root vertices with $N$ smaller than the number of vertices. 
Such $N$-rooted maps turn out to be in bijection with Feynman diagrams for $N$-point functions of a certain quantum field theory. More precisely, as shown in \cite{KriPatVas1},  there is a one-to-one correspondence between $N$-rooted ribbon graphs, or maps, with $e$ edges  and the  $N$-point Feynman diagrams with $(e - N + 1)$ loops   in the so-called {\it scalar quantum electrodynamics}, a quantum field theory which involves quantum fields of two types: a charged scalar field and a photon field. 
Let us emphasize that the connection between Feynman diagrams of this theory and ribbon graphs is of a different nature than in 't  Hooft's work  \cite{thooft} and in matrix models.  
\\

 The established bijection is then used in \cite{KriPatVas1} to obtain explicit expressions and relations for the generating functions of $N$-rooted maps and for the numbers of $N$-rooted maps with a given number of edges without regard to genus using the path integral approach of the quantum field theory. For example, as a generalization of \eqref{ABm1}, the number of $2$-rooted maps of all genera with $e$ edges is given by 
\begin{equation*}
 m_2(e)=\sum_{k=0}^e (-1)^k 
 \sum_{\substack{\mu_1 + \ldots + \mu_{k+1}=e+1 \\ \mu_i \neq 0 } }   \mu_{k+1} ~\prod_{j=1}^{k+1}  \,~(2 \mu_j-1)!!  
-   ~\frac{1}{2} 
  \sum_{k=1}^{e-1} m_1(k) m_1(e-k)\,.
\end{equation*}
Similarly to \eqref{m1Cgn}, the number $m_N(e, v;g)$ of  $N$-rooted graphs of genus $g$ with $v$ vertices and $e$ edges is given by \cite{KriPatVas2}:
\begin{equation}
\label{mN}
m_N(e, v;g) = ~\sum_{\substack{d_1+\dots+d_v=2e \\ d_i \geq 1}} 
 \frac{d_1\cdots d_N}{(v-N)!}\frac{C_{g,v}(d_1,\dots, d_v)}{d_1\cdots d_v}\,.
\end{equation}
In order to obtain the number $m_N(e)$ of $N$-rooted maps with $e$ edges, one can sum the numbers $m_N(e, v;g)$ from \eqref{mN} over $v$ from $N$ to $e+1$ and over $g$ from $0$ to the integer part of $\frac{1+e-v}{2}$. Note that for a given choice of $N,$ the  minimum possible  value $e$ 
may take is $N-1$.\\

 The counting  of rooted maps is of interest in the study of other quantum field theories, some of which having  supersymmetry, see for example \cite{Castro1,Castro2,Castro3,deMello,Prunotto,Vera1,Vera2,Vera3}. 
 \\

In this work we consider $N$-rooted plane trees and study in detail the various subsets of $N$-rooted trees defined by specifying degrees of some of the root vertices.
The numbers of $N$-rooted plane trees with given  degrees of root vertices  have relations to combinatorial objects such as Dyck and lattice paths (see Section \ref{sect_examples}).  In addition, our results provide the first combinatorial interpretations of some sequences listed in the Online Encyclopedia of Integer Sequences \cite{oeis} as well as introduce sequences not listed in the OEIS.
For example, the number of 2-rooted plane trees with $e$ edges is the number of valleys in all the Dyck paths of length $2(e+1)$. 
Note that rooted plane trees (or one-rooted plane trees, that is the case $N=1$) with $e$ edges are dual graphs of the genus zero rooted maps with one vertex, and thus the number of one-rooted plane trees with $e$ edges is the Catalan number $C_e=C_{0,1}(2e).$ 
\\

The maximally rooted plane trees, the {\it tree dicings}, correspond, under the bijection from \cite{KriPatVas1},  to tree level Feynman diagrams of the scalar quantum electrodynamics which determine the dominant approximation of the theory. 
\\

We achieve the enumeration of plane multi-rooted trees by deriving a recursion relation on the numbers of  such trees, which allows us to reduce the enumeration to the ``smaller'' cases, that is trees having fewer roots and smaller degrees of root vertices. This technique is both standard and powerful, going back to Tutte \cite{Tutte0} and Walsh and Lehman \cite{WalshLehman-1, WalshLehman-2}. The same technique is also used in \cite{Dumi, Motohico_Catalan} and leads to the proof  in \cite{Zo, KaZo-Virasoro} of the fact that the generating function of the numbers of dessins d'enfant satisfies Virasoro constraints and the KP hierarchy. 
\\

The paper is organized as follows.  In Section \ref{sect_definitions} we define plane rooted and multi-rooted trees and introduce  notation for the families of trees that we enumerate. In Section \ref{sect_recursions}, in Theorems \ref{thm_recursion} and \ref{thm_simplerecursion}, we derive two recursion relations on the numbers of $N$-rooted trees. The simpler recursion, the one from Theorem \ref{thm_simplerecursion}, is used in Section \ref{sect_formulas} to establish the closed-form expression for the numbers of plane $N$-rooted trees with specified degrees of the root vertices. In Section \ref{sect_examples}, we give a non-exhaustive list of known integer sequences that coincide with the sequences formed by numbers of rooted and multi-rooted trees in various families.  In Section \ref{sect_generating}, we study generating functions of the numbers of $N$-rooted trees with specified degrees of root vertices. In Proposition \ref{prop_Grecursion}, we derive equations allowing us to express such generating functions recursively in terms of those for smaller values of $N$. We give explicit expressions for the cases $N=1,2,3\,.$
Finally, in Section \ref{sect_hyper}, we explore  implications of the recursion from Theorem \ref{thm_recursion} deriving a new binomial identity in Proposition \ref{prop_binomial}  and conjecturing  what we believe to be an original relation between two sums involving binomial coefficients and a corresponding  identity for  certain hypergeometric functions.

\section{$N$-rooted plane trees}

\label{sect_definitions}

A plane tree, that is a tree embedded into a plane, is a special case of a ribbon graph. A ribbon graph, or a fat graph, or a map, is a connected graph with  a fixed cyclic ordering on the set of edges incident to each vertex. More precisely, we have the following definition. 

\begin{definition}
\label{def_tree}
A {\rm  plane tree} is the data $\Gamma = (H, \alpha,\sigma)$ consisting of  a set of half-edges $H = \{h_1,\dots, h_{2e}\}$ with $e$ a positive integer and two permutations $\alpha, \sigma \in S_{2e}$ on the set of half-edges such that
	\begin{itemize}
	\item $\alpha$ is a fixed point free involution,
	\item the subgroup of $S_{2e}$ generated by $\alpha$ and $\sigma$ acts transitively on $H$, 
	\item the number of cycles of $\sigma$ is equal to $1+e$, that is one plus the number of cycles of $\alpha.$
	\end{itemize}
\end{definition}

\begin{figure}[H]
\centering
\includegraphics[scale=0.45]{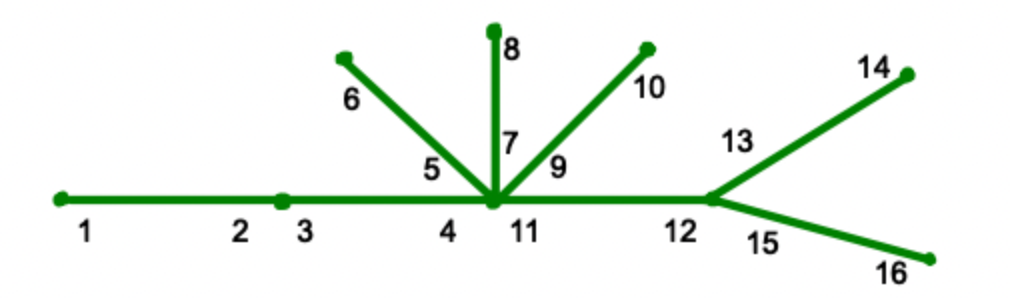}
\caption{The plane tree corresponding to $\alpha=(1\,2)(3\,4)(5\,6)(7\,8)(9\,10)(11\,\,12)(13\,\,14)(15\,\,16)$ 
\\ and $\sigma=(1)(2\,3)(4\,\,11\,\,9\,7\,5)(6)(8)(10)(12\,\,15\,\,13)(14)(16)$}
 \label{fig_tree16}
\end{figure}

\noindent The cycles of  $\alpha$  are transpositions pairing two half-edges to form an edge. Cycles of the permutation $\sigma$ are in bijection with vertices of $\Gamma$, each cycle giving the ordering of half-edges at the corresponding vertex, see Figure \ref{fig_tree16}. 
The transitivity of the group $\langle \sigma, \alpha \rangle$ on the set of half-edges implies the connectedness of the graph $\Gamma$. 
The third condition in Definition \ref{def_tree} ensures that the graph is a tree. 
%

The length of a cycle of $\sigma$, that is the number of half-edges incident to the corresponding vertex, is called the {\it degree} of the  vertex. 
\\

The ordering of the half-edges at every vertex given by $\sigma$ gives the unique way to embed the tree $\Gamma$ into the plane or the sphere. This is done by choosing that each cycle of $\sigma$ corresponds to a vertex such that if one goes around it counterclockwise, then the half-edges attached to this vertex are met in the order given by the cycle of $\sigma$.
\\

Definition \ref{def_tree} without the third condition on the number of cycles of $\sigma$ becomes a definition of a {\it ribbon graph}, or a {\it map}. In this case, some sequence of edges may form a loop and the graph may not necessarily be embedded into a plane without self-intersections.  The cycles of the permutation $\sigma^{-1}\circ \alpha$ correspond to faces of the graph. By gluing a topological disc to each face, we obtain a compact oriented surface into which the graph is embedded. The genus of this surface is called the {\it genus of the ribbon graph}. Thus the genus of a tree is zero. 
\\

In graph enumeration, to avoid double counting, we need to specify which graphs are considered identical. We say that two plane trees are isomorphic if one can be obtained from the other by renumbering the half-edges, and we identify isomorphic trees. We thus obtain trees which may have nontrivial automorphisms, or symmetries. As already mentioned in the introduction, the presence of a distinguished half-edge, called {\it root}, rules out all nontrivial symmetries thus simplifying the task of counting all possible trees. 

\begin{remark}
\label{rmk1}
{\rm
As an exceptional case, a graph consisting of a single vertex and no edges is also considered a plane rooted, or one-rooted, tree.}
\end{remark}

Here we are interested in enumerating $N$-rooted plane trees. Let us start by giving a precise definition. 
\begin{definition}
\label{def_Ntrees}
An $N$-rooted tree is a plane tree, $\Gamma = (H, \alpha, \sigma)$, with the choice of $N$ distinct elements  of $H$, called {\rm  root half-edges}, or {\rm  roots}, belonging to $N$ distinct cycles of $\sigma$, that is incident to $N$ distinct  vertices, called {\rm  root vertices}. The root vertices are labeled by $N$ distinct labels $v_1, \dots, v_N$. The root vertex $v_1$ is distinguished and is called {\rm the first root vertex}; the root at $v_1$ is called {\it the first root} of $\Gamma$.
\end{definition}

\begin{figure}[H]
\centering
\includegraphics[scale=0.5]{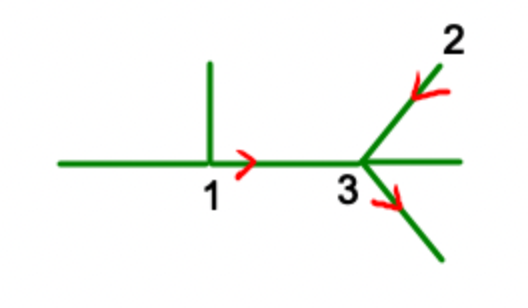}\qquad \qquad
\includegraphics[scale=0.5]{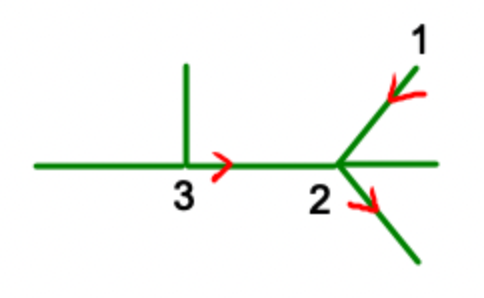}
\caption{Two 3-rooted trees}
 \label{fig_tree}
\end{figure}

\noindent In other words, an $N$-rooted tree is plane tree in which $N$ distinct vertices are chosen and assigned tags $v_1, \dots, v_N\,.$ Moreover, at each of the chosen vertices, an arrow is placed on one of the half-edges incident to the vertex, see Figure \ref{fig_tree}. 

Two $N$-rooted trees isomorphic in the sense of the following definition are identified. 

\begin{definition}
\label{def_tree_iso}
Two plane trees  $\Gamma = (H, \alpha,\sigma)$ and $\Gamma' = (H', \alpha',\sigma')$ are  isomorphic if  there is a bijection $\varphi$ between the sets of half-edges, $\varphi:H\to H',$ such that $\alpha=\varphi^{-1}\alpha'\varphi$ and $\sigma=\varphi^{-1}\sigma'\varphi\,.$ The map $\phi$ is an isomorphism between $\Gamma$ and $\Gamma'\,.$ If both trees  are $N$-rooted, then $\phi$ is an isomorphism of $N$-rooted trees if it maps the $k$th root of $\Gamma$ to the $k$th root of $\Gamma'\,,$ that is if $\phi(r_k)=\phi(r_k')\,.$
\end{definition}

\noindent 
In other words, an isomorphism between two $N$-rooted trees is an isomorphism of plane trees that preserves the labelling of the $N$ root vertices and maps roots to roots.
The only automorphism of an $N$-rooted tree is the identity. For example, the 3-rooted trees from Figure \eqref{fig_tree} are different rooted trees. 
\\

 We denote the number of $N$-rooted plane trees with $e$ edges by $T_N(e)$. In our notation the Catalan numbers are $T_1(e)=C_e$. Let us also denote by $S_N(e)$ the set of all $N$-rooted plane trees with $e$ edges so that $T_N(e)$ is the number of elements in the set $S_N(e)\,.$ Note that due to Remark \ref{rmk1}, we have $T_1(0)=1\,.$
\\

In the case of one rooted trees, one can consider the subset, $S_1(e; d)$, of $S_1(e)$ where the root vertex is specified to be of degree $d$, and the corresponding number of such trees $T_1(e; d).$ By definition, the numbers $T_1(e; d)$ sum up to the $e$th Catalan number $T_1(e)$:
\begin{equation}
\label{dsum}
\sum_{d=1}^e T_1(e; d) = T_1(e)
\end{equation}
and thus we get a well known natural partition of each Catalan number forming the Catalan triangle, see Example \ref{example_triangle} in Section \ref{sect_formulas}. Moreover, for the parts $T_1(e; d)$ of such partitions, we have

\begin{equation}
\label{c1}
T_1(e; d) = \sum_{b = d-1}^{e-1} T_1(e-1; b)
\end{equation}
and 
\begin{equation}
\label{c2}
T_1(e; d)  =\frac{d}{e} \binom{2e-d-1}{e-1}
\end{equation}
for all positive integers $e$ and $d$ with $e\geq d$ and 
with the convention $T_1(0; 0) = 1\,.$ These equations follow from our Theorems \ref{thm_simplerecursion} and  \ref{thm_formula}, respectively. %
It is also possible to verify relations \eqref{dsum} and \eqref{c1} from the closed form expressions \eqref{c2} using the following identities from Section 5.2 of \cite{concrete}:
\begin{equation*}
\sum_{k=1}^{n} k \binom{m-k-1}{m-n-1} = \binom{m}{n-1} \qquad \text{and}\qquad
\sum_{k=0}^n \binom{m-k}{n-k} = \binom{m+1}{n} . 
\end{equation*}
 The sequences of numbers $T_1(e; d)$ for a given $d$ coincide with various other combinatorially interesting sequences, see Section \ref{sect_formulas}. 
 \\

 More generally, one can analogously consider the numbers $T_N(e; d_1, \dotsc, d_k)$ of $N$-rooted plane trees with $e$ edges for which the degrees of $k$ of the $N$ root vertices are fixed to be $d_1, \dotsc, d_k$ respectively, for $0 \leq k \leq N$, as well as the corresponding set of trees $S_N(e; d_1, \dotsc, d_k)$. Here we assume that $d_k$ stands for the degree of the vertex labeled  $v_k$. Note that the number $T_N(e; d_1, \dotsc, d_k)$ does not change if the degrees $d_1, \dots, d_k$ are permuted, for example, $T_9(10; 1,2,3)=T_9(10; 3,1,2)\,.$ The corresponding sets of trees are obtained from one another by relabeling the root vertices. For example, the tree on the left in the Figure \ref{fig_tree} belongs to the set $S_3(6;3,1,4)$ and the tree on the right in the same figure is from the set $S_3(6;1,4,3)$.
\\

In the following lemma, we list relations between degrees $d_i$ and the numbers of edges and roots in a tree. We define the number $T_N(e; d_1, \dotsc, d_k)$ to be zero if the conditions given in this lemma are not satisfied. 
\begin{lemma}
\label{lemma_conditions}
For $N$-rooted trees from the set $S_N(e; d_1, \dotsc, d_k)$ with $k\leq N$, the integer quantities $e, \, N, \, d_i$ satisfy
\begin{eqnarray*}
&&e\geq N-1, \qquad N\geq 1; \\
&&0<d_i\leq e \quad\text{except for the case } N=1, \,\, e=d_1=0 \quad\text{where }\quad T_1(0;0)=1;\\
&&  \sum_{i=1}^{k} d_i \leq e+k-1 \quad\text{where} \quad k=N=e+1 \quad\text{if and only if} \quad \sum_{i=1}^{N} d_i=2e;\\
&& \text{if}\,\,\, e\geq 2\,\,\,\text{then}\,\,\, \Big |\{i\,|\,d_i=1, \, 1\leq i \leq k\}\Big | \leq e\,; \quad\text{if}\,\,\, e=1\,\,\,\text{then}\,\,\,\Big |\{i\,|\,d_i=1, \, 1\leq i \leq 2\}\Big | \leq 2\,.
\end{eqnarray*}
Here the vertical bars denote the number of elements in the set. 
If degrees of  $k-1$  of the $N$ root vertices are fixed to be  $d_1, \dots, d_{k-1}$,  then the highest possible value $D_k$ of the degree $d_k$  of the $k$th root vertex is
\begin{equation}
\label{Dk}
D_k=e+k-1-\sum_{i=1}^{k-1} d_i\,.
\end{equation}
\end{lemma}
{\it Proof.} The condition $e\geq N-1$ holds because the number of vertices in a tree with $e$ edges is $e+1$, so we cannot have more than $e+1$ roots. We set $N\geq 1$ as we do not consider non-rooted trees in this paper. The condition $0<d_i\leq e$ follows again from the fact that our graphs are trees and the rest of the conditions of the second line is the convention stated in Remark \ref{rmk1}. 
\\

If the number of roots is $e+1$, then the tree is maximally rooted, or is a dicing from \cite{WalshLehman-1}. In this case, summing the degrees of all the vertices, we count every edge twice and thus obtain twice the number of edges. If the degrees of only $k\leq N$ root vertices are given, then in the
sum  of  the $k$ given degrees, only the edges connecting two  of the $k$ root vertices are counted twice. As there is at most $k-1$ of such edges, we obtain the condition in the third line. 
\\

The condition in the fourth line follows from the fact that a tree with $e\geq 2$ edges has at most $e$ leaves, that is vertices of degree one. 
\\

Now suppose degrees of $k-1$ of the root vertices are fixed to be $d_1, \dots, d_{k-1}$ and we want to obtain a tree with these degrees and such that the $k$th root vertex be of the maximal possible degree $D_k$. Due to the condition in the third line, 
$\sum_{i=1}^{k-1} d_i+D_k\leq e+k-1$ and thus we want to show that there is a tree for which the equality is attained, that is a tree of $e$ edges for which the degrees $d_1, \dots, d_{k}$ of $k$ root vertices satisfy $\sum_{i=1}^{k-1} d_i+d_k = e+k-1$. To construct such a tree, let us start with $k-1$ disconnected star trees each having a vertex of one of the degrees $d_1, \dots, d_{k-1}$ at its center and all other vertices being of degree one. We then attach all these trees together by glueing together $k-1$ vertices of degree one, one vertex per star tree, thus forming a connected tree and a new vertex of degree $k-1$. This new connected tree has $\sum_{i=1}^{k-1} d_i$ edges. Attaching $e-\sum_{i=1}^{k-1} d_i$ new edges to the vertex of degree $k-1$, we obtain a desired tree containing a vertex of degree $k-1+e-\sum_{i=1}^{k-1} d_i$.
$\Box$
\\

Similarly to the $N=1$ case, one also has the possibility of obtaining $T_N(e; d_1, \dotsc, d_s)$ from $T_N(e; d_1, \dotsc, d_k)$ for $s<k\leq N$, by summing over the appropriate degrees. For example,
for $d_1, \dotsc, d_{N-1}>0$, and $e>0$, as a direct consequence of definition of $T_N(e; d_1, \dotsc, d_k)$,  we have 
\begin{equation*}
T_N(e; d_1, \dotsc, d_{N-1}) = \sum_{d_N = 1}^{D_N} T_N(e; d_1, \dotsc, d_N)\,,
\end{equation*}
where  $D_N$ is the highest possible value of $d_N$ given by \eqref{Dk} with $k=N$. Similarly, summing over more than one $d_i$ to obtain $T_N(e; d_1, \dotsc, d_s)$ for $0 \leq s < N$, we have
\begin{equation}
\label{tns}
T_N(e; d_1, \dotsc, d_{s}) = \sum_{d_{s+1} = 1}^{D_{s+1}} \,\, \sum_{d_{s+2} = 1}^{D_{s+2}}\dots \sum_{d_N = 1}^{D_N} T_N(e; d_1, \dotsc, d_N)\,,
\end{equation}
where $D_k$ is the highest possible value of $d_k$ given by \eqref{Dk}.\\

In the following sections we find closed form expressions for the $T_N(e; d_1, \dotsc, d_s)$. We also derive recursion relations satisfied by the $T_N(e; d_1, \dotsc, d_s)$ and show that our closed form expressions satisfy these relations. \\

\section{Recursion relations for $N$-rooted trees}

\label{sect_recursions}

\noindent The number of $N$-rooted plane trees with a given number of edges can be obtained as recursive combinations of the numbers of rooted trees with fewer edges.  For example, in the case of one rooted maps, the numbers $T_1(e; d)$ can be obtained recursively from the relation:
\begin{equation}
T_1(e; d) = \sum_{k_1 + k_2 = e-1} T_1(k_1; d-1)\, T_1(k_2) = \sum_{k = 0}^{e-1} T_1(k; d-1)\, T_1(e-k-1)\,,
\end{equation}
where $T_1(n)$ are the Catalan numbers. 
More generally, for the numbers $T_N(e; d_1, \dotsc, d_N)$, we prove the following theorem. 

\begin{theorem}
\label{thm_recursion}
For a finite set of non-negative integers $I=\{a_1, \dots, a_l\}$, let $T_N(e; n_1, \dots, n_k, I)$ denote the number $T_N(e; n_1, \dots, n_k, a_1, \dots, a_l)$ with $k$ and $l$ non-negative integers such that $k+l\leq N\,.$
Then, the numbers $T_N(e; d_1, \dotsc, d_N)$  of  $N$-rooted plane trees defined in Section \ref{sect_definitions}  satisfy the recursion
\begin{multline}
\label{recursion}
T_N(e; d_1, \dotsc, d_N) = \sum_{e_1 + e_2 = e-1} \ \biggl( \qquad \sum_{\mathclap{\substack{ I \cup J = \{ d_2, \dotsc, d_N\}\\ I\cap J=\emptyset }}} \qquad T_{|I|+1}(e_1; d_1-1, I)\, T_{|J|+1}(e_2; J)\\ +  \sum_{r=2}^N\qquad\quad\sum_{\mathclap{\substack{I \cup J  = \{d_2, \dotsc, \hat d_r, \dotsc, d_N\} \\ I\cap J=\emptyset}}} \qquad d_r \, T_{|I|+1}(e_1; d_1-1, I)\, T_{|J|+1}(e_2; d_r-1, J) \biggr)\, . 
\end{multline}
Here  $\lvert S\rvert$  stands for the number of elements in the set $S$, 
a hat put over an element of a set signifies that the element is omitted, and
the parameters in each $T_M(e; n_1, \dots, n_k, S)$ need to satisfy conditions of Lemma \ref{lemma_conditions} for the $T_N$ to be nonzero.   
When summing over all partitions of the set of degrees into two disjoint sets $I$ and $J$, the degrees $d_2, \dots, d_N$ are considered as labels of the root vertices and not as integers, so that even if $d_i=d_j$, a partition for which $d_i\in I$ and $d_j\in J$ is different from a partition for which $d_j\in I$ and $d_i\in J\,.$

\end{theorem}

{\it Proof.} We prove this recursion by establishing a bijection between the set $S_N(e; d_1, \dotsc, d_N)$ of $N$-rooted trees counted by the number $T_N(e; d_1, \dotsc, d_N)$  and the set of trees counted by the right hand side of \eqref{recursion}. To this end, consider a tree from the set $S_N(e; d_1, \dotsc, d_N)$ and denote its root vertices by $v_1, \dots, v_N$ so that $v_j$ is of degree $d_j$. We put this tree in correspondence with a set of trees with one fewer edges by removing the edge containing the root half-edge at $v_1$. This separates the tree into two trees disconnected from each other. A similar strategy was used to prove analogous recursions, for example, in \cite{WalshLehman-1} and  \cite{Motohico_Catalan}. There are two essentially different cases: when the edge in question connects $v_1$ to another root vertex and when it connects $v_1$ to a non-root vertex. These two cases correspond to the two terms in the right hand side of recursion \eqref{recursion}.\\

{\bf Case 1.} Let $\Gamma$ be a tree from the set $S_N(e; d_1, \dotsc, d_N)$ such that its edge containing the root edge at $v_1$ connects $v_1$ to a non-root vertex. Denote this edge by $l$. We remove $l$ from $\Gamma$ and put an arrow on the half-edge following $l$ at $v_1$ in the counterclockwise order (if such a half-edge exists), thus creating a new root at $v_1$. At the same time, we put an arrow on the half-edge following $l$ in the counterclockwise order at the other vertex incident to $l$, thus creating a new root vertex denoted by $v$. Thus two rooted trees are created, $\Gamma_1$ with $e_1$ edges containing the root vertex $v_1$ of degree $d_1-1$ and $\Gamma_2$ with $e_2$ edges containing the root vertex $v$ of unknown degree that can take any possible value; here $e_1+e_2=e-1$. For a given split of the set of degrees  $\{d_2, \dots, d_N\}$ into two disjoint sets $I$ and $J$ according to the vertices included in $\Gamma_1$ and those included in $\Gamma_2\,,$ the number of such pairs of trees $(\Gamma_1, \Gamma_2)$ is thus $T_{|I|+1}(e_1; d_1-1, I)\, T_{|J|+1}(e_2; J)$. Note that the set $S_{|J|+1}(e_2; J)$ contains all $(|J|+1)$-rooted trees with degrees given by the set $J$ while the degree of the $(|J|+1)$th vertex is not prescribed. This latter vertex is our vertex $v$.  Summing over all possible splits of the set of root vertices into two parts, we obtain all possible pairs of trees   $(\Gamma_1, \Gamma_2)$ that can be created in this case, as well as the first line in \eqref{recursion}.      \\

{\bf Case 2.} Let $\Gamma$ be a tree from the set $S_N(e; d_1, \dotsc, d_N)$ such that its edge containing the root at $v_1$ connects $v_1$ to another root vertex $v_r$.  In this case, we first remove the arrow at $v_r$ and then apply the procedure from Case 1 with $v=v_r$. We thus create two trees $\Gamma_1$ with $e_1$ edges containing the root vertex $v_1$ of degree $d_1-1$ and $\Gamma_2$ with $e_2$ edges containing the root vertex $v_r$ of degree $d_r-1\,,$ with $e_1+e_2=e-1$.  For a fixed splitting of the set $\{d_2, \dots, d_N\}\setminus\{d_r\}$ into two disjoint sets $I$ and $J$, the number of created pairs of trees is $d_r\,T_{|I|+1}(e_1; d_1-1, I)\, T_{|J|+1}(e_2; d_r-1, J)$.  The factor of $d_r$ is due to the fact that we started by removing the arrow at $v_r$ and thus the resulting trees are the same no mater which of the $d_r$ half-edges at $v_r$ carried an arrow in $\Gamma$. Again, summing over all possible disjoint sets $I$ and $J$ and then over all values of $r$ between $2$ and $N$, we obtain the second line of  \eqref{recursion}. \\

Reciprocally, let $N_1,\,N_2$ be two positive integers such that $N_1+N_2=N+1$ or $N_1+N_2=N$ and let $e_1, \;e_2$ be two non-negative integers such that $e_1+e_2+1=e$. Starting with an $N_1$-rooted tree $\Gamma_1$ with $e_1$ edges and with a $N_2$-rooted tree $\Gamma_2$  with $e_2$ edges,  we can join them by a new edge between their respective first root vertices and create a tree from the set $S_N(e; d_1, \dotsc, d_N)$. In this case, the new edge must be inserted  so that it follows the two roots in the clockwise order at the respective root vertices and the first root of $\Gamma_1$ must be moved to the new half-edge at the same vertex. As for the first root of $\Gamma_2$, it must be removed if $N_1+N_2=N+1$. If  $N_1+N_2=N$, the first root vertex of $\Gamma_2$ must be renamed, its root removed and then placed at each half edge at the same vertex, including the new half-edge, thus creating a set of $d+1$ distinct graphs, where $d$ is the degree of the first root vertex of $\Gamma_2\,.$ In the case with $N_1+N_2=N+1$, we create a tree from the set $T_N(e; d_1, \dotsc, d_N)$ whose first root belongs to an edge connecting a root vertex to a non-root vertex, whereas in the case $N_1+N_2=N$, the first root belongs to an edge which connects two root vertices of the created tree. 
$\Box$
\\

Recursion \eqref{recursion} has a structure similar to that of the recursion for the generalized Catalan numbers $C_{g,v}(d_1, \dots, d_v)$ obtained in \cite{WalshLehman-1} and rederived in \cite{Motohico_Catalan}. Since in \cite{Motohico_Catalan} the recursion for $C_{g,v}(d_1, \dots, d_v)$ was also linked to the topological recursion of Chekhov-Eynard-Orantin from \cite{EO}, a natural question to ask is whether our recursion \eqref{recursion} can also be obtained by the topological recursion applied to some algebraic curve. 
\\

However, it is difficult to derive closed form expressions for the number of $N$-rooted trees from the recursion of Theorem \ref{thm_recursion}. Another, simpler, recursion given in the next theorem turns out to be more useful in this respect.

\begin{theorem}
\label{thm_simplerecursion}
The numbers $T_N(e, d_1; \dotsc, d_s)$ \eqref{tns} of  $N$-rooted plane trees with $e$ edges, where  $s$ of the $N$ ordered root vertices are specified to be of degrees $d_1, \dotsc, d_s$, respectively, with $d_i \geq 0$, satisfy the following recursion
\begin{multline}
\label{simplerecursion}
T_N(e; d_1, \dotsc, d_{N-1}, d_N) = \sum_{i=1}^{N-1} d_i T_{N-1}(e-1;d_1, \dots, d_{i-1}, d_i+d_N-2,d_{i+1}, \dots, d_{N-1}) 
\\
+ T_N(e-1; d_1, \dots, d_{N-1}) - \sum_{d=1}^{d_N-2}T_N(e-1; d_1, \dots, d_{N-1}, d)
\,. 
\end{multline}
Here we assume that $T_M(n; n_1, \dots, n_k)$ is zero if the conditions of Lemma \ref{lemma_conditions} are not satisfied. 
\end{theorem}

{\it Proof.} This recursion is obtained in a way similar to the proof of Theorem \ref{thm_recursion}. This time, a bijection between the set $S_N(e; d_1, \dots, d_N)$ of $N$-rooted trees counted by the number $T_N(e; d_1, \dots, d_N)$ and the set of trees counted by the right hand side of  \eqref{simplerecursion} is established by contracting an edge instead of removing it.  Let $\Gamma$ be a tree from the set $S_N(e; d_1, \dots, d_N)$ and denote its root vertices by $v_1, \dots, v_N$. The edge we contract is the one containing the root half-edge at the $N$th root vertex $v_N$ of  $\Gamma$.  Let us denote this edge by $l$. There are again two cases.
\\

{\bf Case 1. } The edge $l$ connects $v_N$ to a non-root vertex $v$. Let us move the arrow of the root half-edge at $v_N$ to the half-edge following $l$ in counterclockwise order at $v_N$ if $d_N\neq 1$ and at $v$ if $d_N= 1$ and then contract the edge $l$. We thus create a new root vertex by merging $v_N$ and $v$. Its degree depends on the degree of $v$ and thus can be any number between $d_N-1$ and the maximal degree possible given the constraint of the other degrees $d_1, \dots, d_{N-1}$ and the total number of edges $e-1$ of the new tree. Performing this procedure for every tree of the set $S_N(e; d_1, \dots, d_N)$, we thus obtain all trees of the set $S_N(e-1; d_1, \dots, d_{N-1})$ where the degree of the $N$th root vertex is not specified except for the trees where this degree is smaller than $d_N-1\,.$ Thus, we get in total the number of trees given by the second line in \eqref{simplerecursion}.
\\

{\bf Case 2. } The edge $l$ connects $v_N$ to another root vertex $v_i\,.$ We first remove the arrow from the root half-edge at $v_i$ and then apply the procedure of Case 1. Let us call the vertex obtained by merging $v_N$ and $v_i$ again by $v_i\,.$ Its degree is $d_i+d_N-2\,.$ Thus this procedure applied to all trees of the set $S_N(e; d_1, \dots, d_N)$ satisfying assumption of Case 2, yields the number of trees counted by the first sum in the right hand side of \eqref{simplerecursion}. The factor of $d_i$ takes care of the fact that we started by removing the arrow at $v_i$, and thus there are $d_i$ different trees that produce identical resulting tree. 
\\

Note that the arrows help us to keep track of where the contracted edge used to be so that we can reverse the procedure similarly to the proof of Theorem \ref{thm_recursion}. 
$\Box$

\section{Counting $N$-rooted trees}

\label{sect_formulas}

In this section we give the closed form expressions for the numbers $T_N(e; d_1, \dotsc, d_N)$ and show that they satisfy the recursion relation in Theorem \ref{thm_recursion}. 
\begin{theorem}
\label{thm_formula}
Under the conditions of Lemma \ref{lemma_conditions}, 
the numbers  $T_N(e; d_1, \dotsc, d_N)$ of  $N$-rooted plane trees with $e$ edges and degrees $d_1, \dotsc, d_N$ of the  $N$ ordered root vertices
 are given by
\begin{equation}
\label{formula}
T_N(e; d_1, \dotsc, d_N) = \frac{(e-1)!}{(e+1-N)!} \binom{2e - 1- \sum_{i=1}^N d_i}{e-N} \prod_{j=1}^N d_j \,\qquad \text{if}\quad e\neq 0 
\end{equation}
and the only non-zero case with $e=0$ being
\begin{equation*}
T_1(0;0)=1\,.
\end{equation*}
\end{theorem}
For the case of $N=1$, these formulas reduce to the numbers given by \eqref{c2}. On the other hand, in \cite{WalshLehman-1}, the following result was obtained for the number of genus $g$ maximally rooted maps (dicings) having
one face and degree $d_j$ of the vertex $v_j$: 
\begin{equation*}
F_g(d_1, \dots, d_v) = \frac{(v+2g-2)!}{2^{2g}} \left( \prod_{j=1}^v d_j\right) \sum_{\substack{k_1+\cdots + k_v= g \\ k_1, \ldots, k_v\geq  0}} \,\,\prod_{j=1}^v\frac{1}{2k_j+1}\binom{d_j-1}{2k_j}\,.
\end{equation*}
In the case of genus zero, these numbers count maximally rooted trees and give 

\begin{equation*}
F_0(d_1, \dots, d_v) = (v-2)! \left( \prod_{j=1}^v d_j\right) \binom{-1}{0}^v \,.
\end{equation*}
Binomial coefficients with a negative upper entry  are calculated using the identity 
\begin{equation*}
\binom{-\alpha}{\beta} = (-1)^\beta \binom{\alpha+\beta -1}{\beta} 
\end{equation*} which is valid whenever $\alpha$ and $\beta$ are non negative integers. This leads to $\binom{-1}{0} = \binom{0}{0} = 1$. Therefore
the result for $F_0(d_1, \dots, d_v)$  coincides with  \eqref{formula} after identifying $v=N=e+1$ and  $d_1+\dots+d_N=2(N-1)$, and using the convention $\binom{-1}{-1}=1\,.$ \\

{\it Proof of Theorem \ref{thm_formula}.}
In order to prove that the numbers $T_N$ are given by \eqref{formula}, we prove that expressions \eqref{formula} satisfy the recursion from Theorem \ref{thm_simplerecursion}.  This recursion allows to construct all the numbers $T_N$ starting from the base case of one rooted tree with one edge and degree $1$ of the root vertex, that is the number $T_1(1;1)$. For this base case, formula \eqref{formula} is valid as it gives $T_1(1;1)=1\,,$ note that $\binom{0}{0}=1\,.$ Thus it remains to prove that expressions in the right hand side of \eqref{formula} satisfy \eqref{simplerecursion}, which we rewrite using \eqref{tns} in the form
\begin{equation*}
T_N(e; d_1, \dots, d_N) \hspace{-0.1cm}=\hspace{-0.2cm}\sum_{i=1}^{N-1}\hspace{-0.1cm} d_i T_{N-1}(e-1;d_1, \dots, d_{i-1}, d_i+d_N-2, d_{i+1},\dots, d_{N-1}) 
+\hspace{-0.3cm} \sum_{d=d_N-1}^{D_N} \hspace{-0.3cm} T_N(e-1; d_1, \dots, d_{N-1}, d)\,.
\end{equation*}
Here $D_N=e+N-2-\sum_{j=1}^{N-1} d_j$ is the maximal possible value of $d$ as given in \eqref{tns} with $e$ replaced by $e-1\,.$ Plugging in \eqref{formula} in the right hand side of this equality, we obtain
\begin{multline}
\label{temp} 
T_N(e; d_1, \dots, d_N) =\frac{(e-2)!}{(e-N+1)!}\binom{2e-1-\sum_{j=1}^{N}d_j}{e-N}      \sum_{i=1}^{N-1}(d_i+d_N-2)  \prod_{j=1}^{N-1}d_j  
\\
+\frac{(e-2)!}{(e-N)!}   \sum_{d=d_N-1}^{D_N}d\binom{2e-3-d-\sum_{j=1}^{N-1}d_j}{e-N-1}  \prod_{j=1}^{N-1}d_j  .
\end{multline}
Our objective now is to show that the right hand side coincides with the right hand side of \eqref{formula}. Rewriting the last sum using the following constant
\begin{equation}
\label{change}
m=2e-2-\sum_{j=1}^{N-1}d_j,
\end{equation}  
we obtain
\begin{equation*}
\sum_{d=d_N-1}^{D_N}d\binom{2e-3-d-\sum_{j=1}^{N-1}d_j}{e-N-1} =
\sum_{d=d_N-1}^{D_N}d\binom{m-d-1}{m-D_N-1}\,.
\end{equation*}
Introducing a new summation variable $k=m-d$, we have
\begin{multline*}
\sum_{d=d_N-1}^{D_N}d\binom{2e-3-d-\sum_{j=1}^{N-1}d_j}{e-N-1}
=\sum_{k=m-D_N}^{m-d_N+1}(m-k)\binom{k-1}{m-D_N-1}\\
=m\sum_{k=m-D_N}^{m-d_N+1}\binom{k-1}{m-D_N-1}- \sum_{k=m-D_N}^{m-d_N+1}k\binom{k-1}{m-D_N-1}\\
=m\sum_{k=m-D_N-1}^{m-d_N}\binom{k}{m-D_N-1}- (m-D_N)\sum_{k=m-D_N}^{m-d_N+1}\binom{k}{m-D_N}\,,
\end{multline*}
where in the last line, in the first sum,  we changed the summation variable from $k$ to $k-1$, and in the second sum, we used the  property $\binom {a+1}{b+1} =\frac{a+1}{b+1} \binom {a}{b}$ of binomial coefficients.  Applying the following property 
of  Pascal's triangle
\begin{equation*}
\sum_{k=a}^{b}\binom{k}{a}=\binom{b+1}{a+1} 
\end{equation*}
which can be found, for example, in \cite{concrete}, Table 174, we obtain 
\begin{multline*}
\sum_{d=d_N-1}^{D_N}d\binom{2e-3-d-\sum_{j=1}^{N-1}d_j}{e-N-1} =
m\binom{m-d_N+1}{m-D_N}- (m-D_N)\binom{m-d_N+2}{m-D_N+1}
\\
=  \frac{(m-D_N)(d_N-1)+D_N}{m-D_N+1}\binom{m-d_N+1}{m-D_N}\,,
\end{multline*}
where we used again the property $\binom {a+1}{b+1} =\frac{a+1}{b+1} \binom {a}{b}$ in the second term. Going back to the original notation \eqref{change} and using the last result in \eqref{temp}, we see that the binomial coefficients in both terms of \eqref{temp} coincide. By pulling out the common factor and a straightforward calculation, we prove that the right hand side of \eqref{temp} coincides with the right hand side of \eqref{formula}, which finishes the proof. 
$\Box$

\begin{corollary}
\label{cor_rec}
The numbers $T_N(e; d_1, \dots, d_N)$ given by  Theorem \ref{thm_formula} satisfy the recursion in Theorem \ref{thm_recursion}.
\end{corollary}

{\it Proof.} This follows from the fact that the numbers of trees $T_N(e; d_1, \dots, d_N)$ satisfy both recursions of Theorems \ref{thm_recursion} and \ref{thm_simplerecursion}. 
$\Box$

\begin{remark}
\label{remark_2N}
As can easily be checked from \eqref{formula}, for $N\leq e-1$, the following numbers  of $N$-rooted trees with degrees of all $N$ root vertices specified coincide: $T_N(e; 2, \dots, 2,1)=T_N(e; 2, \dots, 2)\,.$ The condition $N\leq e-1$ is necessary to ensure that both sets of parameters satisfy conditions of Lemma \ref{lemma_conditions}. 

\end{remark}

\begin{proposition}
\label{prop_sum} Let $0\leq s\leq N$ and $e>0$ be integers and the conditions of Lemma \ref{lemma_conditions} be satisfied.
The sums of  $T_N(e; d_1, \dotsc, d_N)$  over the values of $d_{s+1}, \dots, d_N$ defined in \eqref{tns} are given by
\begin{equation}
\label{sum}
T_N(e; d_1, \dotsc, d_s) = \frac{(e-1)!}{(e+1-N)!} \binom{2e +N-s-1- \sum_{j=1}^{s} d_j}{e+N-2s} \prod_{j=1}^s d_j \ .
\end{equation}
\end{proposition} 
{\it Proof.} We prove this by induction on the number $N-s$ of the degrees $d_j$ which we sum over. Note that when this number is zero, that is when $s=N$, formula \eqref{sum} coincides with \eqref{formula}, giving us the base case of the induction. Now, assume that the statement of the proposition is true for some value $r\leq N$ taking the place of $s$ in \eqref{sum}. We then want to prove \eqref{formula} for $s=r-1\,.$ By definition \eqref{tns}, 
\begin{equation*}
T_N(e; d_1, \dotsc, d_{r-1}) = \sum_{d_r = 1}^{D_r} T_N(e; d_1, \dotsc, d_r)\,,
\end{equation*}
where $D_r=e+r-1-\sum_{i=1}^{r-1} d_i$ is the maximal possible value of $d_r\,.$
By the induction hypothesis, the $T_N$ in the sum are given by \eqref{sum} with $s=r$ and thus
\begin{equation*}
T_N(e; d_1, \dotsc, d_{r-1}) =\frac{(e-1)!}{(e+1-N)!}  \prod_{j=1}^{r-1} d_j  \sum_{d_r = 1}^{D_r} d_r\binom{2e +N-r-1- \sum_{j=1}^{r} d_j}{e+N-2r}   \,.
\end{equation*}
It remains to use the following binomial identity, see \cite{concrete}, p. 177: 
\begin{equation*}
\sum_{k = 1}^{p} k\binom{m-k-1}{m-p-1}=\binom{m}{p-1} =\binom{m}{m-p+1}    
\end{equation*}
with $m=2e+N-r-\sum_{j=1}^{r-1} d_j$, $k=d_r$ and $p=D_r\,.$
$\Box$\\

As an immediate corollary of Proposition \ref{prop_sum}, we obtain the number of $N$-rooted trees with a given number of edges. 
\begin{corollary} Let $0\leq N-1\leq e\,.$ The number of $N$-rooted trees with $e$ edges is given by 
\begin{equation}
\label{tn}
T_N(e) =  \frac{(e-1)!}{(e+1-N)!} \binom{2e  + N- 1}{e +N} =   \frac{(e-1)!}{(e+1-N)!} \binom{2e  + N- 1}{e - 1} \,, \quad\text{if}\quad e>0
\end{equation}
with the only nonzero case for $e=0$ being 
\begin{equation*}
T_1(0) =1\,.
\end{equation*}
\end{corollary}
{\it Proof.}
Summing the numbers $T_N(e; d_1, \dotsc, d_{N})$ given by \eqref{formula} over all the degrees $d_1, \dotsc, d_N$, that is setting $s=0$ in Proposition \ref{prop_sum},  we obtain expression \eqref{tn} for $T_N(e)$. 
$\Box$

Setting $N=1$ in \eqref{tn} we get
\begin{equation}
T_1(e) = \frac1{e} \binom{2e}{e-1} = 
\frac1{e+1} \binom{2e}{e} = C_e\,,
\end{equation}
which is just the expression for Catalan numbers. \\

\section{Some explicit examples of our sequences and relations to known sequences}
\label{sect_examples}

In this section, we give a few samples of our sequences and mention alternative combinatorial interpretations for some of them. 

\begin{example}
{\rm  The following sequences  of numbers of $N$-rooted trees can be obtained from \eqref{tn}.}

\medskip

\begin{tabular}{|l|c | c|c|c|c|c|c|c|c|c|clc|c|}
 \hline
 \diagbox{$N$}{$e$} & 0 & 1 & 2  &3 &4 &  5 &6 & 7  &8 &9 \\
 \hline
 1 &1 & 1 &2&5& 14& 42& 132  &429 &1 430  & 4 862\\
 \hline
 2 & & 1  &  5 & 21 & 84  & 330 & 1 287 &  5 005 & 19 448 &  75 582\\
 \hline
  3 & - & -   &  6 & 56 & 360  &  1 980 & 10 010 & 48 048 & 222 768 &  1 007 760 \\
 \hline
   4 & - & -   &  - & 72 & 990   &  8 580  &60 060   & 371 280   & 2 116 296  &11 395 440  \\
 \hline
  5 & - & -   &  - & - &  1 320  &  24 024  & 262 080 & 2 227 680  & 16 279 200  & 107 442 720 \\
 \hline
 6  & - & -   &  - & - &  -   &  32 760 & 742 560 &  9 767 520  & 97 675 200  & 823 727 520  \\
 \hline
  7  & - & -   &  - & - &  -   &- &1 028 160   &  27 907 200 &    429 770 880 &   4 942 365 120 \\
 \hline
 \end{tabular}
\end{example}

\bigskip

\begin{figure}[H]
\centering
\includegraphics[scale=0.35]{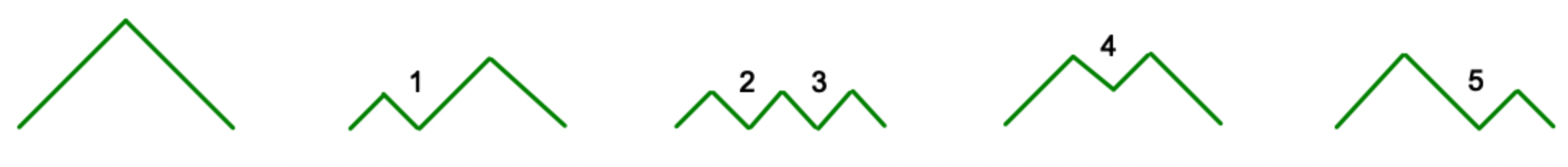}
\caption{Five valleys in Dyck paths of length 6}
 \label{fig_Dyck}
\end{figure}

The sequence $T_2(e)$ in the second line coincides with the OEIS \cite{oeis} sequence \href{https://oeis.org/A002054}{(A002054)}, which has several combinatorial descriptions, one of which being the total number of valleys in all Dyck paths of length $2(e+1)$.  For example, as shown in Figure \ref{fig_Dyck} , there is a total of 5 valleys in all  Dyck paths of length 6, which corresponds to $T_2(2)=5.$
If we divide the sequence $T_3(e)$ by $3$, we obtain the OEIS sequence \href{https://oeis.org/A074922}{(A074922)}. 
 The sequences for $N >3$ are not recorded in OEIS.  
The first non-zero entries in each line of the above table, in other words the values $T_N(N-1)$ form the sequence \href{https://oeis.org/A001763}{(A001763)}. In our interpretation, this is a sequence of numbers of maximally rooted trees, or the numbers of trees among dicings. 
\\

For a given $N$, by fixing values of degrees $d_j$ in the numbers $T_N(e; d_1, \dotsc, d_s)$, we obtain various sequences indexed by $e\in \mathbb N\,.$ Most of these sequences are new, but some  are known in other  contexts. 
For example,  the numbers $T_1(e; 1)$ and $T_1(e; 2)$ give the Catalan numbers $C_{e-1}$. The numbers $T_1(e; 3)$ for varying $e$ give the OEIS sequence \href{https://oeis.org/A000245}{(A000245)}, $T_1(e; 4)$ is the OEIS sequence \href{https://oeis.org/A002057}{(A002057)}, $T_1(e; 5)$ is the OEIS sequence \href{https://oeis.org/A00344}{(A00344)}, $T_1(e; 6)$ is the OEIS sequence \href{https://oeis.org/A003517}{(A003517)}. 
\\

\begin{example}
\label{example_triangle}
{\rm  Given that several sequences $T_1(e; d)$ for fixed values of $d$  coincide with various known sequences, we list here some of them. One recognizes the Catalan triangle, which is therefore now given a new interpretation in terms of one rooted tree graphs. Once the first two rows (corresponding to $d=0$ and $d=1$) are given, all other entries can be determined recursively using 
$T_1(e,d) = T_1(e,d-1) - T_1(e-1,d-2)$.} 

\medskip

\begin{tabular}{|l|c| c|c|c|c|c|c|c|c|c|c|c|}
 \hline
 \diagbox{$d$}{$e$} & 0 & 1 & 2  &3 &4 &  5 &6 & 7  &8 &9 & 10     & 11 \\
 \hline
 0 & 1 & 0 & 0 & 0 & 0 &0 & 0 & 0 & 0 & 0  &0       & 0 \\
 \hline
 1 &-    & 1 &  1  & 2  & 5 &  14& 42& 132  &429 &1 430 & 4 862      &   16 796  \\
 \hline
 2 &   - & -  & 1&  2  & 5  & 14  & 42 & 132 &  429 & 1 430 &  4 862      &   16 796 \\
 \hline
  3 & - & -   &  -  & 1  & 3  & 9  & 28 &  90 & 297 &   1 001  & 3 432       & 11 934  \\
 \hline
   4 & - & -   &  - & -  &1  & 4 & 14 & 48 & 165 &  572  & 2 002       &  7 072 \\
 \hline
  5 & - & -   &  - & - &  - & 1  & 5& 20 & 75 &  275 & 1 001       & 3 640   \\
 \hline
 6  & - & -   &  - & - &  -   &  - & 1 & 6 &27 & 110 &  429      & 1 638  \\
 \hline
  7  & - & -   &  - & - &  -   &- &- & 1 & 7 & 35  &  154   & 637 \\
 \hline
 \end{tabular}
\end{example}

\bigskip

Our sequences $T_1(e;d)$  for $d \geq 3 $ have  interesting connections with certain classes of  North-East  lattice  paths ({\it i.e.} paths with steps $(0,1)$ or $(1,0)$)  going from $(0,0)$ to $(n,n)$. 
For example $T_1(e;3)$ is  given by the number of such paths with $n=e-1$ which  bounce off  the diagonal to the right ({\it i.e.} touch the diagonal after going up and then go to the right) only once but never cross the diagonal vertically \cite{PanRemmel}. 
\\

$T_1(e;4) $  is given by the number of paths in a square with sides of length $e-1$ which bounce exactly twice to the right off the diagonal but never cross it vertically, and so on. So for $d \geq 3$, $T_1(e;d)$  is given by the number of paths in a square with sides of length $e-1$ which bounce exactly $d-2$ to the right off the diagonal but never cross the diagonal vertically. In the last section, we will show how the generating function of these paths given in \cite{PanRemmel} can be used to construct the generating function for the $T_1(e;d)$, including those with $d=1$ and $d=2$.  
\\

Many other interpretations can be given to the values $T_1(e;d)$ in terms of lattice paths. As a last example, the numbers $T_1(e;4)$ is also equal to the number of such paths in squares of sides $e-2$ that have exactly one horizontal crossing of the diagonal and no vertical crossings  \cite{PanRemmel}.

\begin{example}
{\rm  A sample of the sequences $T_2(e;d_1,d_2)$ for $d_1 \leq d_2$ (note that $T_2(e;d_1,d_2)$ are symmetric under the exchange of $d_1$ and $d_2$). Note that the sequences in lines 2 and 3 coincide for $e\geq 3$ as in Remark \ref{remark_2N}.}

\medskip

\begin{tabular}{|l|c | c|c|c|c|c|c|c|c|c|clc|c|}
 \hline
 \diagbox{$d_1, d_2$}{$e$}  & 1 & 2  &3 &4 &  5 &6 & 7  &8 &9 & 10\\
 \hline
1, 1 & 1 &1 & 3 & 10&   35 &  126 &  462 & 1 716 & 6 435 & 24 310 \\
 \hline
 1, 2  &-    &  2& 4& 12& 40 & 140 & 504 & 1 848 & 6 864 & 25 740    \\
 \hline
 2, 2 &   - & -  & 4 & 12& 40 & 140 & 504 & 1 848 & 6 864 & 25 740    \\
 \hline
  1, 3  & - & -   &  3 & 9 &  30 &  105 & 378 & 1 386 & 5 148 &  19 305   \\
 \hline
   2, 3  & - & -   &  - & 6  & 24 & 90 & 336 & 1 260 &  4 752  & 18 018  \\
 \hline
3, 3 & - & -   &  - & - &  9 & 45 & 189 & 756 & 2 970 &  11 583  \\
 \hline
 1, 4 & - & -   &  - & 4 &16 & 60 & 224 &  840 & 3 168 &   12 012\\
 \hline
2, 4 & - & -   &  - & - &  8 &40  &  168  & 672 &  2 640 & 10 296   \\
 \hline
 \end{tabular}
\end{example}

The sequence   $T_2(e; 1, 1)$ is the OEIS sequence \href{https://oeis.org/A088218}{(A088218)} (or the closely related \href{https://oeis.org/A001700}{(A001700)}) which gives the total number of leaves in all rooted ordered trees with $n$ edges.  

The sequence  $T_2(e; 1, 2)$ corresponds to the twice central binomial coefficients, OEIS sequence \href{https://oeis.org/A028329}{(A028329)}, which has several 
combinatorial interpretations, for example as the number of North-East lattice paths from $(0,0)$ to $(n+1,n+1)$ that cross the diagonal an even number of times\cite{PanRemmel}. 
\\

The rows of the table correspond to various sequences that can be found on the OEIS website,  after dividing by the product $d_1 d_2$ all the terms. For example, the sequence $T_2(2,3)$ is six times the sequence \href{https://oeis.org/A001791}{(A001791)}.

\begin{example}
\label{example_degreeone}
{\rm  A sample of the sequences $T_N(e;d_1,d_2,\ldots, d_N)$   with all degrees equal to one, $d_1=\ldots=d_N=1$.}

\medskip

\begin{tabular}{|l|c | c|c|c|c|c|c|c|c|c|clc|c|}
 \hline
 \diagbox{$N$}{$e$}  & 1 & 2  &3 &4 &  5 &6 & 7  &8 &9 & 10\\
 \hline
1 & 1 &1 & 2 & 5 & 14 & 42 & 132 &   429 & 1 430 &  4 862\\
 \hline
 2  &1 &    1 &  3 &  10 & 35 &  126 & 462 & 1 716 & 6 435 & 24 310  \\
 \hline
 3 &   - & -  & 2 & 12 &  60 &   280 & 1 260 & 5 544 & 24 024 &102 960   \\
 \hline
 4 & - & -   &-   & 6 & 60 & 420 &  2 520 &  13 860 & 72 072 & 360 360  \\
 \hline
 5  & - & -   &  - & -   & 24 & 360 & 3 360 & 25 200 & 166 320 & 1 009 008   \\
 \hline
6 & - & -   &  - & - &  - & 120 & 2 520 &  30 240 & 277 200 & 2 162 160   \\
 \hline
 7 & - & -   &  - & - &- & -  & 720 & 20 160 & 302 400 &  3 326 400 \\
 \hline
8  & - & -   &  - & - & -& - & -& 5040 &   181 440 & 3 326 400 \\
 \hline
 \end{tabular}
\end{example}

The numbers $T_N(e;d_1, \ldots , d_N) $ with all the degrees equal to $1$ are of interest. The sequence $T_1(e;1)$ is of course made of the Catalan numbers and $T_2(e;1,1)$ has already been discussed above. After dividing by $(N-1)!$, the sequence $T_3(e;1,1,1)$ , $T_4(e;1,1,1,1)$  and $T_5(e,1, \ldots,1)$ can be found in the OEIS. In the case of $T_5$, our result provides the first combinatorial interpretation of the sequence. 
\\

 The sequences $T_N(e;1,\ldots ,1)$  after division by $(N-1)!$,  are given  given by the coefficients $e(q,m)$ presented in \cite{Cossali} (see the columns of their Table 1), which were introduced through a recursion formula and they also appear in the triangle 
 \href{https://oeis.org/A088617}{(A088617)}.
 \\

For even $N$, after dividing by $(N-1)!$, the sequences $T_N(e;1, \ldots , 1)$ also  appear as the columns in the triangle given in the OEIS entry \href{https://oeis.org/A281000}{(A281000)} and give the first combinatorial interpretation of these numbers.
\\

The sequences $T_3(e;d_1,d_2,d_3)$ have entries in the OEIS only for low values of the degrees. 
\\

Not only do our numbers of  rooted trees provide new integer sequences, they also unify a large number of sequences found in the OEIS that are given  extremely disparate descriptions.

 \section{Generating functions for $N$-rooted plane trees}

\label{sect_generating}

Finding an expression for a generating function of a given sequence is sometimes easier than determining closed form expressions for the terms of the sequence. The approach to map enumeration through generating functions is used very often, see for example \cite{AB, BC, BCR, KriPatVas1, JV, WalshLehman-1, WalshLehman-2, KaZo-Rational} and many other works. The expression for the generating function of the Catalan numbers was used in \cite{Motohico_Catalan} to inform the choice of the algebraic curve for the Eynard-Orantin topological recursion which can be used to produce the generalized Catalan numbers $G_{g,v}(d_1, \dots, d_v)$ defined in the introduction.
\\

Let us first consider  the generating function of the numbers of one-rooted trees $T_1(e)$, which is the same as the generating function of the Catalan numbers $C_e\,.$ Let us denote this function by $C(t)$ and define it as follows
\begin{equation}
\label{Ct-def}
C(t) = \sum_{e=1}^\infty  T_1(e) \, t^e = \sum_{e=1}^\infty  C_e \, t^e.
\end{equation}
This function satisfies 
\begin{equation}
\label{Ct-rel}
C(t) = 1+t\,C^2(t)
\end{equation}
as can be seen using the interpretation of the coefficients in \eqref{Ct-def} as numbers of rooted trees similarly to the proof of Theorem \ref{thm_recursion}. Namely, the set $S_1$ of all one-rooted trees can be produced from two copies of itself by taking one tree from each copy of the set, removing the two roots and creating a new edge between the two former root vertices such that the new edge precedes the former root edges in the counterclockwise order at both vertices. In fact, this statement needs a slight correction - the one-rooted tree with one vertex and no edges cannot be produced in this way. In the new edge, the half-edge on the side of the graph from the ``first'' copy of the set $S_1$ becomes the root of the obtained tree.  Translating this process into the terms of the generating function, we obtain \eqref{Ct-rel}, where creating a new edge corresponds to increasing the exponent of the variable $t$ in the series \eqref{Ct-def}, and thus to multiplying the generating function by $t$. 
\\

Having obtained \eqref{Ct-rel}, we can now solve this equation for $C(t)$ and obtain the well known form of the generating function of the Catalan numbers:
\begin{equation}
\label{Cat} 
C(t) =  \frac{1-\sqrt{1-4t}}{2t}\,.
\end{equation} 
Note that the other solution, the function $\hat C(t) =  \frac{1+\sqrt{1-4t}}{2t}$ contains the same information as it is related to \eqref{Cat} by $C(t)\hat C(t) = 1/t\,.$
Let us now define a generating function $G_1(t,x)$ of the numbers $T_1(e;d)$ by
\begin{equation}
\label{G1-def}
 G_1(t,x) =   \sum_{d=0}^\infty  \sum_{e=d}^\infty  T_1(e;d) \, x^d t^e.
\end{equation}
Applying the same logic as in \eqref{Ct-rel}, one can see that $G_1(t,x)$ satisfies
\begin{equation}
\label{G1-rel} 
 G_1(t,x) = 1+t \, x \,G_1(t,x) G_1(t,1)\,.
\end{equation} 
Namely, the set $S_1$ without the degenerated tree that has no edges can be obtained from two copies of $S_1$ exactly as in \eqref{Ct-rel}. The difference is that now we keep track of the degree of the first root vertex and thus by creating a new edge we also increase the degree of the root vertex, thus multiplication by $x$ on the right in \eqref{G1-rel}. Moreover, the root of the tree from the ``second'' copy of $S_1$ disappears and thus the degree of its root vertex is unimportant, which corresponds to setting $d=1$ for the second factor of $G_1$ on the right in \eqref{G1-rel}. 

Note now that $G_1(t,1) = C(t)$ and thus we can solve \eqref{G1-rel} for $ G_1(t,x)$ and obtain
\begin{equation}
\label{G1} 
 G_1(t,x) = \frac{1}{1-t \, x\,C(t)} = \frac{1 +\sqrt{1-4t}}{1+\sqrt{1-4t}-2x  \, t }\,.
\end{equation} 
In general, for $1\leq s \leq N$, we can introduce the following generating functions
\begin{equation}
\label{GN-def}
 G_N(t,x_1, \dots, x_s) =   \sum_{d_1=0}^\infty \dots \sum_{d_s=0}^\infty \,  \sum_{e=0}^\infty  T_N(e;d_1, \dots, d_s) \, x_1^{d_1}\cdots x_s^{d_s} t^e,
\end{equation}
where  $T_N(e;d_1, \dots, d_s)$ vanishes unless the parameters satisfy the conditions of Lemma \eqref{lemma_conditions}. Note that $ G_N(t,x_1, \dots, x_s)$ and  $G_N(t,x_1, \dots, x_N)$ with $1\leq s \leq N$ are related by  $G_N(t,x_1, \dots, x_s) = G_N(t,x_1, \dots, x_s, 1, \dots, 1)$ where the variables $x_{s+1}, \dots, x_{N}$ are set to 1.  It thus suffices to know $G_N(t,x_1, \dots, x_N)$ to know all the functions in \eqref{GN-def}. Due to the symmetry of the numbers $T_N(e;d_1, \dots, d_s)$ under permutation of the degrees $d_1, \dots, d_s$, the functions $ G_N(t,x_1, \dots, x_s)$ are symmetric in the variables $x_1, \dots, x_s\,.$
\\

A generalization of \eqref{G1-rel}, allows us to find expressions for all generating functions $G_N(t, x_1,\dots,  x_N)$ recursively from those with smaller values of $N$. Namely, we have the following recursion. 
\begin{proposition}
\label{prop_Grecursion}
The  generating functions $G_K(t, x_1,\dots, x_s)\,$   defined by \eqref{GN-def}  satisfy
\begin{eqnarray}
\label{diffG}
&G_N(t,x_1, \dots, x_{N-1}) &= \left( 2t\frac{\partial}{\partial t} - \sum_{k=1}^{N-1}x_k \frac{\partial}{\partial x_k}\right) G_{N-1}(t,x_1, \dots, x_{N-1})\,;
\\
\nonumber
&G_N(t; x_1, \dotsc, x_N) &= tx_1  \quad \sum_{\mathclap{\substack{ I \cup J = \{ x_2, \dotsc, x_N\}\\ I\cap J=\emptyset }}} \qquad G_{|I|+1}(t; x_1, I)\, G_{|J|+1}(t; J)\\ 
\label{Grecursion}
&\qquad\qquad\qquad\qquad+& tx_1\sum_{r=2}^N\qquad\quad \sum_{\mathclap{\substack{I \cup J  = \{x_2, \dotsc, \hat x_r, \dotsc, x_N\} \\ I\cap J=\emptyset}}} \qquad x_r \, G_{|I|+1}(t; x_1, I)\, \frac{\partial}{\partial x_r}\biggl\{x_r\,G_{|J|+1}(t; x_r, J) \biggr\}\, . 
\end{eqnarray}
The sum in the second equation is taken over all possible splits of the given sets into two disjoint sets $I$ and $J$; a hat put over an element of a set signifies that the element is omitted.  
\end{proposition}
{\it Proof.} Let us start by proving the first equation of the proposition. Recall that $G_N(t,x_1, \dots, x_{N-1})=G_N(t,x_1, \dots, x_{N-1}, 1)$ is the generating function of the numbers $T_N(e; d_1, \dots, d_{N-1})$ where we keep track only of the degrees of the first $N-1$ root vertices of the $N$-rooted trees with $e$ edges.  We can obtain the trees of the corresponding set $S_N(e; d_1, \dots, d_{N-1})$ from the trees of the set $S_{N-1}(e; d_1, \dots, d_{N-1})$ by introducing the $N$th root in all possible ways. There are $2e$ half-edges in total and $d_1+\dots, +d_{N-1}$ of the half-edges cannot be chosen to be a new root. Thus there are $2e-d_1-\ldots -d_{N-1}$ ways to choose the $N$th root in each tree of the set $S_{N-1}(e; d_1, \dots, d_{N-1})$ and each of these choices produces a distinct tree in the set $S_N(e; d_1, \dots, d_{N-1})$ since the presence of one root eliminates all non-trivial automorphisms of the graph. In other words,  $$T_N(e; d_1, \dots, d_{N-1}) = (2e-d_1-\ldots -d_{N-1})T_{N-1}(e; d_1, \dots, d_{N-1}).$$ Thus if in the series defining $G_{N-1}(t,x_1, \dots, x_{N-1})$ we multiply every term containing $t^ex_1^{d_1}\cdots x_{N-1}^{d_{N-1}}$ by $2e-d_1-\ldots-d_{N-1}$, we obtain the series for $G_N(t,x_1, \dots, x_{N-1}).$ This is exactly what is done by applying the differential operator in the right hand side of \eqref{diffG} to $G_{N-1}(t,x_1, \dots, x_{N-1})$.
\\

The second equation of the proposition is nothing but a rewriting of the recursion of Theorem \ref{thm_recursion} in terms of the generating functions. To see the validity of \eqref{Grecursion}, it is best to perform the inverse of the procedure described in the proof of Theorem \ref{thm_recursion}. Namely, we start with two rooted trees, call one of them ``left'' and another one ``right'' and join them by a new edge connecting the first root vertices of the two trees. The new edge is added in such a way that it precedes the root half-edge of both trees in the counterclockwise order at both vertices. Then the arrows marking the first root half-edges are removed from both trees and the half-edge of the new edge that is incident to the vertex of the tree ``on the left'' is chosen for the new first root. 
In this way, we create a new rooted tree. Case 1 of the proof of Theorem \ref{thm_recursion} gives a new tree with the first root  connecting a root vertex to a non root vertex while Case 2 gives a new tree with the first root connecting two root vertices. These two cases correspond to the two terms in the right hand side of \eqref{Grecursion}. The factor of $t$ which increases the exponent of $t$ in every term in the series defining the generating function by one, corresponds to adding a new edge. The factor of $x_1$ corresponds to the increasing of $d_1$ due to the added edge. In Case 2, we also put an arrow at the former first root vertex of the tree ``on the right'' in all possible ways. The number of such ways is one plus the degree of the former first root vertex of the tree ``on the right''. Multiplying of the terms of the generating function by such factor is achieved by the operator $\frac{\partial}{\partial x_r}x_r$ where $x_r$ is the variable corresponding to the root vertex in question. Summing over all disjoint splittings of the set $\{ x_2, \dotsc, x_N\}$, we include all possible distributions of the root vertices among the two trees that we connect by a new edge. 
$\Box$
\\

As an example, here are the equations of Proposition \ref{prop_Grecursion} allowing to find $G_2(t, x_1, x_2)$ from $G_1(t,x)$: 
\begin{eqnarray*}
&&G_2(t,x) = \left( 2t\frac{\partial}{\partial t} - x \frac{\partial}{\partial x}\right) G_1(t,x)\,;
\\
&&G_2(t,x_1,x_2) = x_1 t\left( \,  G_2(t,x_1,x_2) \, G_1(t,1) +   G_1(t,x_1)\,  G_2(t,x_2) +   x_2  \, G_1(t,x_1)  \frac{ \partial}{\partial x_2} \biggl\{ x_2 \,  G_1(t,x_2) \biggr\} \right)\,.
\end{eqnarray*}
\\
Starting from expression \eqref{G1} for the generating function $G_1(t,x)$ and using the first of the above equations, we obtain $G_2(t,x)=G_2(t,x,1)\,.$ This gives us all the necessary ingredients to obtain $G_2(t, x_1, x_2)$ from the second of the above equations.   This leads to
\begin{equation}
\label{G2}
G_2(t,x_1,x_2) = \frac{   t \,  x_1  \, x_2 \,  (y+1)^5  }{ 2 y \, (1+ y -2t x_1)^2 \, (1+ y -2t x_2)^2}\,\,,
\end{equation}
where $y$ stands  for the square root appearing in the generating function \eqref{Cat} of the Catalan numbers:
\begin{equation*}
y=\sqrt{1-4t}, \qquad 1-y^2=4t\,.
\end{equation*}

Having found $G_2$, we apply the equations of Proposition \ref{prop_Grecursion} with $N=3$, to obtain the following expression for $G_3:$
 \begin{multline*}
 G_3(t,x_1,x_2,x_3) = \frac{ t^2 \, x_1 \, x_2 \, x_3 \, (y+1)^{10}~}{2 \, y^3 \,  (1+y-2t x_1)^3 \, (1+y-2tx_2)^3 \, (1+y-2t x_3)^3}  \Biggl(
1-y + (2y-2t) (x_1+x_2+x_3)  \\ + \frac{1}{2} (y-1)(1+3y-6t) (x_1 x_2 + x_1 x_3 + x_2 x_3) 
 + \frac{1}{4} (y-1)^2 (2+4y-10t) x_1 x_2 x_3 \Biggr) .
 \end{multline*} 
 This leads us to the following conjecture. 
\begin{conjecture}
The generating functions defined by \eqref{GN-def} for $N\geq 3$ have the form
 \begin{eqnarray*}
 G_N(t,x_1,\ldots , x_N) &=& \frac{(N-1)!}{4} ~ \frac{t^{N-1} \,  (1+y)^{N^2+1} \prod_{j=1}^N x_j}{y^{2N-3} \, \prod_{k=1}^N (1+y-2 t \, x_k)^N}~  {\cal{P}}_N (y,x_1, \ldots x_N) \,,
 \end{eqnarray*}
where ${\cal{P}}_N (y,x_1, \ldots x_N)$ is a polynomial in all its variables. 
\end{conjecture}

We have noticed that the generating function  $F_{2,4}(x_2,x_4,t)$ introduced in  \cite{PanRemmel} with the second argument set to zero contains 
all the $T_1(e,d)$ after some manipulations. Interestingly, that function generates the number of north-east lattice paths, according to the number of times the path bounces to the right off the diagonal (organized in powers of $x_2$) and the number of times it crosses the diagonal vertically (counted by the powers of $x_4$). This function is given by
\begin{eqnarray}
 F_{2,4} \bigl(t,x_2,x_4\bigr) = \frac{(x_2-2) (-1 +\sqrt{1-4 t} ) + 2 (x_2-1) t}{x_4 (-1+\sqrt{1-4t}) + (2+ x_2 (-1 + \sqrt{1-4t}) +3x_4 - x_4 \sqrt{1-4t}   ) t } .
 \end{eqnarray}
Note that this  corrects  a typo in \cite{PanRemmel}  (the factor of $x_2$ in the denominator was misplaced).  With $x_4=0$, this generating function reproduces all the $T_1(e,d) $ for $ d \geq 3$, which implies that the one-rooted trees count certain types of lattice paths (see Section 5 for more details). For $d=1$ and $d=2$, some rearrangement must be made to obtain our generating function $G_1(t,x)$ \eqref{G1-def}, \eqref{G1} from $F_{2,4}$. More precisely, we have the following relation
\begin{eqnarray*} G_1 (t,x) &=& x^2 t    \bigl( F_{2,4}(t,x_2=x,x_4=0) +1 - 2   C(t)  \bigr) + x^2 t \bigl( C(t)-1 \bigr)  + x \, t \,  C(t)  + 1  \,.
\end{eqnarray*}
\\

In the table of Example \ref{example_degreeone}, we presented some of the sequences $T_N(e;d_1, \ldots , d_N)$ with all degrees equal to 1. The corresponding generating functions, let us call them $G_{N}^{1 \ldots 1} (t)$,  can be obtained by differentiating the  generating functions $G_N(t,x_1, \ldots , x_N)$ with respect to all the $x_i$ once and then setting $x_i=0$. For $G_{1}^{1 \ldots 1}(t)$ we of course have the generating function of the Catalan numbers. For $N=2$ and $N=3$, the generating functions are
  \begin{eqnarray*}
 G_{2}^{1 \ldots 1} (t) &=& \frac{1}{2} \frac{(1+\sqrt{1-4t})t}{\sqrt{1-4t}} ,\\
 G_{3}^{1 \ldots 1} (t) &=& \frac{2 t^3}{(1-4t)^{\frac{3}{2}} }.
  \end{eqnarray*}

  \section{Binomial and hypergeometric identities from $T_N(e, I_N)$}
\label{sect_hyper}

In Section \ref{sect_formulas}  we obtained the expression for $T_N(e, d_1, \ldots , d_N)$ by solving the recursion formula of Theorem \ref{thm_simplerecursion}. 
The recursion formula of Theorem \ref{thm_recursion} is much more difficult to work with. Even for a fixed value of $N$, showing that  numbers \eqref{formula} satisfy the recursion of Theorem \ref{thm_recursion}  is 
nontrivial for $N>1$. For $N=2$ the recursion can be proven using formulas found in \cite{concrete}. Starting with $N=3$, new identities for sums of product of binomial coefficients are needed. These sums can be expressed in terms of generalized hypergeometric functions with certain arguments, and therefore the recursion formula can be used to derive new hypergeometric identities. We will illustrate this for $N=3$.
\\

For $N=3$, Theorem \ref{thm_recursion} gives 
\begin{eqnarray}
\label{recu2}
T_3(e;d_1,d_2,d_3) &= &\sum_{e_1+e_2=e-1} \Biggl[ T_1(e_1;d_1-1)  \, T_3(e_2;d_2,d_3) + T_3(e_1;d_1-1,d_2,d_3) \, T_1(e_2) 
\nonumber  \\
&~& + T_2(e_1;d_1-1,d_2) \, T_2(e_2;d_3)  ~+ ~T_2(e_1;d_1-1,d_3) \, T_2(e_2;d_2)   \nonumber   \\
&~& + d_2 \,  T_1(e_1;d_1-1) \,  T_2(e_2;d_2-1,d_3)  ~+ ~d_2 \, T_2(e_1;d_1-1,d_3)  \, T_1(e_2;d_2-1)   \nonumber  \\
 &+&~d_3  \, T_2(e_1;d_1-1,d_2) \, T_1(e_2;d_3-1) ~
+ ~d_3\,  T_1(e_1,d_1-1)  \, T_2(e_2,d_3-1,d_2) \Biggr]\,.\qquad
\end{eqnarray}

If $d_1=1$, the proof of \eqref{recu2} is trivial, the only terms on the right hand side of \eqref{recu2} being nonzero are the first, the fifth (at the condition that $d_2 >1$)  and the last  (at the condition that $d_3 >1$) ones. These three terms may contribute  because  $T_1(0;0)=1$;  their sum  can be checked to be equal to $T_3(e;1,d_2,d_3)$.
\\

 Let us focus on the first sum on the right hand side of (\ref{recu2}), which we will denote $\Sigma_1$  (a function of $e,d_1,d_2 $ and $d_3$).  If we choose to sum over $e_1$, it will  range  from $(e_1)_{min}$ to $ e-1-(e_2)_{min}$. The minimum value of $e_1$ is determined by the conditions of Lemma 1 applied to the factor $T_1(e_1;d_1-1)$ which, taking into account the fact that we are assuming $d_1 >1$,   gives $e_1 \geq d_1-1$. 
 \\
 
  The minimum value of  $e_2$ is determined by the factor $T_3(e_2;d_2,d_3)$ and  depends on the values of $d_2$ and $d_3$. If $d_2=d_3=1$, the condition is  $e_2 \geq 2$  and   if $d_2+d_3 >2$, it is 
  $e_2 \geq d_2+d_3-1$. Let us consider first the case $d_2 + d_3 >2$.
Setting $e_2=e-e_1-1$, the first term of  (\ref{recu2})  becomes, after plugging in expressions \eqref{formula} for the numbers $T_1$ and $T_3$,
\begin{eqnarray}
\Sigma_1 = \sum_{e_1=d_1-1}^{e-d_2-d_3} \frac{(d_1-1)(e-2-e_1)d_2 d_3}{e_1} \binom{2e_1-d_1}{e_1-1} \binom{2e -2e_1-2-d_2-d_3}{e-2-e_1} .
\label{pik}
\end{eqnarray}

If $d_2=d_3=1$, the upper limit on the sum over $e_1$ in \eqref{pik} is equal to  $e-2$. However, note that  the expression evaluated at  $e_1=e-2$ vanishes. Therefore when $d_2=d_3=1$, the sum in \eqref{pik} is effectively up to $e-3$. But this coincides with $ e-1-(e_2)_{min}$ as $(e_2)_{min}=2$ in this case.  Therefore  \eqref{pik} is also valid for the case $d_2=d_3=1\,.$
\\

After making the appropriate changes of variables, one finds 
$$ \Sigma_1 = (d_1-1)(e-2) d_2 d_3 \, S_1(n_1,r_1,s_1) - (d_1-1) d_2 d_3 \, S_5(n_1,r_1,s_1), $$ where $S_1$ and $S_5$ refer to the sums listed in the appendix and 
 $$ n_1 =e+1-\sum_{i=1}^3 d_i \,, ~~~~r_1=d_1-2\,, ~~~~~~~s_1=d_2+d_3-2\,.$$
 Following similar steps, the next seven terms in the right hand side of \eqref{recu2} are found to be 
\begin{eqnarray*}  \Sigma_2  &=&  (d_1-1)(e-2) d_2 d_3 S_4(n_2,r_2) - (d_1-1) d_2 d_3 S_5(n_2,r_2,s_2) + T_3(e-1;d_1-1,d_2,d_3), \\
 \Sigma_3  &=& (d_1-1)  d_2 d_3 S_5(n_3,r_3,s_3),
 \\
  \Sigma_4  &=& (d_1-1) d_2 d_3\, S_5(n_4,r_4,s_4),
  \\ \Sigma_5  &=&  d_2 d_3 (d_1-1)(d_2-1) S_1(n_5,r_5,s_5) ,
  \\ \Sigma_6  &=& d_2 d_3 (d_1-1)(d_2-1) S_1(n_6,r_6,s_6),
  \\ \Sigma_7  &=& d_2 d_3 (d_1-1)(d_3-1) S_1(n_7,r_7,s_7) ,
  \\ \Sigma_8   &=& d_2 d_3 (d_1-1)(d_3-1)  \, S_1(n_5,r_5,s_5) ,
  \end{eqnarray*}
 where
 \begin{eqnarray*}
 n_2 &=&n_1  , ~~~~~~~~~~~~r_2= d_1+d_2+d_3 - 6 , \quad  ~s_2=2,\\
 n_3 &= & n_1 , ~~~~~~~~~~~~r_3=d_1+d_2-4,\quad \quad  \quad  ~~s_3=d_3, \\
 n_4 &= &n_1, ~~~~~~~~~~~~r_4=d_1+d_3-4, \quad \quad  \quad  ~~s_4=d_2, \\
 n_5 &=&n_1+1,~~~~~~~r_5=d_1-2,\quad \quad  \quad \quad  \quad  ~~ \,s_5=d_2+d_3-4, \\
 n_6 &= & n_1+1,~~~~~~~r_6=d_2-2,\quad \quad  \quad \quad  \quad  ~~ s_6=d_1+d_3-4, \\
 n_7 &= & n_1+1,~~~~~~~r_7=d_3-2,\quad \quad  \quad \quad  \quad  ~~ s_7=d_2+d_1-4.
  \end{eqnarray*}

 Adding up all eight terms, $\Sigma_1 $ to $\Sigma_8$,  we rewrite (\ref{recu2})  in the following form, valid under the assumption $d_1 \geq 1$: 
 \begin{eqnarray}
T_3(e;d_1,d_2,d_3) &=& d_1 d_2 d_3  (e-1)  \binom{2e-1- \sum_{i=1}^3 d_i }{e-3} +  (d_1-1) d_2 d_3 \Bigl(  S_5(n_3,r_3,s_3) 
\nonumber
\\&~& ~~~~~~~~~~~~+ S_5(n_4,r_4,s_4)  -\, S_5(n_1,r_1,s_1)
   -\, S_5(n_2,r_2,s_2)
 \Bigr)    .
  \end{eqnarray}
 \\

 From \eqref{formula}, we know that this is equal to $T_3(e;d_1,d_2,d_3) = d_1 d_2 d_3  (e-1)  \binom{2e-1-\sum_{i=1}^3 d_i }{e-3}$ and we   have therefore proven the following proposition. 
 \begin{proposition}
 \label{prop_binomial}
 The quantity $S_5 (n,r,s)=\sum_{k=0}^{n} \binom{2n-2k+s}{n-k} \, \binom{2k+r}{k}$ with $n,r,s$ being arbitrary non-negative integers satisfies

 \begin{eqnarray}
   S_5(n_3,r_3,s_3) + S_5(n_4,r_4,s_4)
 - \, S_5(n_1,r_1,s_1) 
 - \, S_5(n_2,r_2,s_2) =0 . \label{bbb}
 \end {eqnarray}
 \end{proposition}

  From the definition of the $S_5$ as a sum, it  is not obvious at first sight  that \eqref{bbb} is satisfied because the various $S_5$ are evaluated with different arguments.   We have therefore
  obtained a nontrivial combinatorial identity from the recursion formula  \eqref{recu2}. \\
  
  We have been unable to find a closed form formula for the sum $S_5$ in the literature or to derive one. However, 
  we have the following conjecture which  expresses $S_5$ as a different sum and which makes (\ref{bbb}) automatically satisfied. 
 
  \begin{conjecture}
  \label{conj_binom}
  Let $n,r,s$ be  integers and $ n \geq 0$.  Then
 \begin{eqnarray}
   \label{hyper}
  \sum_{k=0}^{n } \binom{2n-2k+s}{n-k} \, \binom{2k+r}{k}  =\sum_{k=0}^{[\frac{n}{2} ] } \binom{2n+2+r+s}{n-2k} ,
 \end{eqnarray}
  where $[\frac{n}{2}]$ denotes the integer part of $\frac{n}{2}$.
  \end{conjecture}
  
  The usefulness of the representation in the right hand side of \eqref{hyper} is that  it depends on $r$ and $s$ through their  sum only, which is not obvious from the expression on the left.
   If this conjecture is  correct, it automatically ensures that (\ref{bbb}) is satisfied since the sums $r_1+s_1$, $r_2+s_2$, $r_3+s_3$ and $r_4+s_4$ are all equal to $\sum_{i=1}^3 d_i-4$. \\

Conjecture \ref{conj_binom} can be expressed as a conjecture relating certain generalized  hypergeometric functions evaluated at $z=1$. The left hand side of \eqref{hyper} can  be shown to be equal to 
  \begin{eqnarray} \binom{2n+s}{n} ~\,_4F_3 \left(1+ \frac{r}{2}, \frac{1+r}{2},-n,-n-s;1+r,-n - \frac{s}{2},-n + \frac{1-s}{2};1\right),    \label{conjec1}
   \end{eqnarray}
   when $r,s \geq 0$. 
   \\

   For the right  hand side of \eqref{hyper}, 
   consider first  the case of odd $n$. Defining $p=n-2k$ , the sum becomes 
 \begin{eqnarray}\sum_{p=1}^n{\vphantom{\sum}}'
  \binom{2n+2+r+s}{p} ,  \label{tutu}  \end{eqnarray} where the prime indicates that the sum is only over the odd values of  $p$ from $1$ to $n$. Using now the identity 
$$ \binom{2n+2+r+s}{p} = \binom{2n+1+r+s}{p}+ \binom{2n+1+r+s}{p-1},$$ we
obtain that the sum \eqref{tutu} can   be written as
$$ \sum_{l=0}^n \binom{2n+1+r+s}{ l }, $$  where now $l$ takes all integer values from $0$ to $n$. 
 One can check that the same result is valid \color{black} when $n$ is even.
  Therefore, we have 
 \begin{eqnarray}\sum_{k=0}^{[\frac{n}{2}]} \binom{2n+2+r+s}{n-2k} = \sum_{l=0}^n \binom{2n+1+r+s}{l}.  \label{yaya}  \end{eqnarray}
 
 \begin{theorem}
 For a non-negative integer $n$ and two integers $r $ and $s$ satisfying $r+s \geq 0$, 
 
 \begin{eqnarray}
  \sum_{l=0}^n \binom{2n+1+r+s}{l} = \binom{2n+2+r+s}{n} ~\,_3F_2 \left(1,\frac{1-n}{2}, -\frac{n}{2}; \frac{n+3+r+s}{2},2+ \frac{n+r+s}{2};1\right).  \label{conjec2}   
 \end{eqnarray}
 
 \end{theorem}

We had conjectured this result and it  has been proven by Professor Robert S. Maier \cite{Maier}. 
The non trivial proof uses  a recursion formula for a certain class of hypergeometric functions which  can be solved  in terms of  sums, one of which is \eqref{yaya}.
\\

  Our conjecture \ref{conj_binom}  can therefore be stated as the equality of  the two hypergeometric functions of \eqref{conjec1} and  \eqref{conjec2}. Note again that this would imply  that \eqref{conjec1}  depends on $r$ and $s$ through their sum only, which is not obvious from the expression. 
 \\
 
 To summarize this section,  the recursion of Theorem \ref{thm_recursion} applied to expressions \eqref{formula} for  $T_3$ implies \eqref{bbb}, which  is a nontrivial condition on the sum $S_5$.  We have been led to conjecture a different form for the sum $S_5$ (the right hand side of \eqref{hyper}), which satisfies \eqref{bbb}. Our conjecture can also be expressed as an equality between two hypergeometric functions. We believe that using Theorem \ref{thm_recursion} applied to  expressions \eqref{formula} for $T_N$ with $N>3$ will lead to ever more complex identities on various sums of products of binomial coefficients,  some  of which should be new results. 

\appendix
\markboth{}{}
\renewcommand{\thesection}{\Alph{section}}
\numberwithin{equation}{section}
\section{Appendix}
\label{appendix_hypergeometric}

Consider the following sums:
\begin{eqnarray*}
S_1 (n,r,s) &: =& \sum_{k=0}^n \frac{1}{k+r+1} \binom{2n-2k+s}{n-k} \, \binom{2k+r}{k} ,\\
S_2 (n,r,s) &:=& \sum_{k=0}^{n+s} \frac{1}{k+r+1} \binom{2n-2k+s}{n-k} \, \binom{2k+r}{k} ,\\
S_3 (n,r)&:=& \sum_{k=0}^{n} \frac{1}{n-k+1} \binom{2n-2k}{n-k} \, \binom{2k+r}{k} ,\\
S_4(n,r) &:=& \sum_{k=0}^n \frac{1}{n-k+1} \binom{2n-2k+2}{n-k} \, \binom{2k+r}{k}, \\
S_5 (n,r,s)&:=& \sum_{k=0}^{n} \binom{2n-2k+s}{n-k} \, \binom{2k+r}{k} ,\\
S_6(n,r,s,t) &:=& \sum_{k=0}^n \frac{r}{tk + r } \binom{tk +r}{k} \, \binom{tn-tk+s}{n-k} .
\end{eqnarray*}
Not all these sums appear in Section \ref{sect_hyper}, but all are closely related to sums we needed and we include them for completeness. 
\\

The sum $S_6$ is given in Equation (5.62) in \cite{concrete}. It is understood that if the parameters considered are such that the factor $tk+r$ is equal to zero for some value of $k$, then in that term the binomial factor $\binom{tk+r}{k}$ is taken to cancel the factor $\frac{1}{tk+r}$. With this convention, $S_6$ is well defined for  $n \geq 0$ and   $r,s,t \in \mathbb{Z}$ and is equal to 
$$ S_6(n,r,s,t) = \binom{tn+r+s}{n} . $$

The sums $S_1$ to $S_4$ can be obtained from $S_6$ after appropriate changes of variables, with the following results:
\begin{eqnarray*}
S_1 (n,r,s) &=&   \frac{1}{r+1} \binom{2n+r+s+1}{n}, \\
S_2(n,r,s) &=& \frac{1}{r+1} \binom{2n+r+s+1}{n+s}, \\
S_3(n,r) &=&  \binom{2n+r+1}{n}, \\
S_4(n,r) &=&  \binom{2n+r+2}{n},
\end{eqnarray*}
where, in $S_1$ and $S_2$,  $n$ and $r$ are integers satisfying $n,r \geq 0$ and $s \in \mathbb{Z}$ (the results are actually valid for a wider range of values but the general results are not needed for this work). The conditions for $S_3 $ and $S_4$ are $ n\geq 0$ and $r \in \mathbb{Z}$.  
\\

We could not find the sum $S_5$  in  closed form in the literature, as mentioned in the main text,  see \eqref{hyper}.
\\

{\bf Acknowledgements.}   We are grateful to Professor Robert S. Maier for several illuminating exchanges on properties of  generalized  hypergeometric functions and techniques to prove related identities.  AG and VS gratefully acknowledge
support from the Natural Sciences and Engineering Research Council of Canada through a Discovery grant and from the Universit\'e de Sherbrooke. GK is thankful to Guido Carlet and to the Universit\'e de Bourgogne, IMB Dijon and IPaDEGAN for making an academic visit possible and to the Universit\'e de Sherbrooke where part of this work was done.   We thank the anonymous referee for  comments that helped us to improve our text and for pointing out important references.

\end{document}